\documentclass{amsart}
\usepackage{graphicx}
\newtheorem{theorem}{Theorem}
\newtheorem{lemma}[theorem]{Lemma}
\newtheorem{corollary}[theorem]{Corollary}
\newtheorem{definition}{Definition}
\newtheorem{proposition}{Proposition}
 
\newcommand{\bla}{\bullet}
\newcommand{\whi}{\circ}

\newcommand{\basic}{\mathcal{B}}

\title{Crossing and alignments of permutations}

\author{Sylvie Corteel}
\address{CNRS PRISM, Universit\'e de Versailles, France 
and Institute Mittag Leffler, Sweden}
\date{May 2, 2005}
\begin{document}
\maketitle

\begin{abstract}
We derive the continued fraction form of the generating
function of some new $q$-analogs of the Eulerian numbers $\hat{E}_{k,n}(q)$ introduced 
by Lauren Williams  building on work of Alexander Postnikov.
They are related to the number
of alignments and weak exceedances of permutations. 
We show how these numbers are related to crossing
and generalized patterns of permutations 
We generalize to the case of decorated
permutations. Finally we show how these numbers appear naturally
in the stationary distribution of the ASEP model. 
\end{abstract} 
\section{Introduction}

The purpose of this paper is to derive the continued fraction form of generating
function of some new $q$-analogs of Eulerian numbers
$\hat{E}_{k,n}(q)$ introduced 
by Lauren Williams \cite{wi} building on work of Alexander Postnikov \cite{po} and to link them
to generalized patterns.
These numbers are related to the number
of alignments and weak exceedances
of permutations and come from the enumeration
of totally positive Grassmann cells. We refer to \cite{po} and \cite{wi}
for details. We start by some definitions.

Let $\sigma=(\sigma(1),\ldots ,\sigma(n))$ be a permutation
of $[n]$. Then let
\begin{itemize}
\item $A_+(i)=\{j\ |\ j<i\le \sigma(i)<\sigma(j)\}$; $A_+(\sigma)=\sum_{i=1}^n |A_+(i)|$.
\item $A_-(i)=\{j\ | \ j>i>\sigma(i)>\sigma(j)\}$; $A_-(\sigma)=\sum_{i=1}^n |A_-(i)|$.
\item $A_{+,-}(i)=
\{j\ |\ j\le \sigma(j)<\sigma(i)<i\}\cup \{j\ |\ \sigma(i)<i<j\le \sigma(j)\}$; $A_{+,-}(\sigma)=\sum_{i=1}^n |A_{+,-}(i)|$.
\end{itemize}
For example let $\sigma=(4,7,3,6,2,1,5)$ then $A_+(\sigma)=3$, $A_-(\sigma)=1$, 
and $A_{+,-}(\sigma)=2$.\\

We can also define these parameters using the permutation diagram. If
$\sigma$ is a permutation of $[n]$. We draw a line and
put the numbers from 1 to $n$ and we draw an edge from
$i$ to $\sigma(i)$ above the line if $i\le \sigma(i)$ and 
under the line otherwise. Then $A_+(\sigma)$ is the number of pairs of
nested edges above (resp. under) the line and $A_{+,-}(\sigma)$
is the number of pairs of edges such that one is above the line and
the other under and such that their support do not intersect.  

\begin{figure}[ht!]
  \begin{center}
    \includegraphics[scale=0.5]{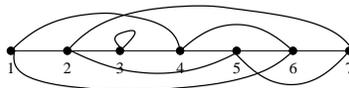}
  \end{center}\centering
  \caption{The permutation diagram of  $\sigma=(4,7,3,6,2,1,5)$}
  \label{perm}
\end{figure}

The permutation diagram of  $\sigma=(4,7,3,6,2,1,5)$ is on Figure \ref{perm}. 
The pairs of edges
contributing to  $A_{+}(\sigma)$ are $\{(1,4),(3,3)\}$, $\{(2,7),(3,3)\}$ and $\{(2,7),(4,6)\}$.
The pair of edges
contributing to  $A_{-}(\sigma)$ is $\{(5,2),(6,1)\}$.
The pairs of edges
contributing to  $A_{+,-}(\sigma)$ are $\{(7,5),(1,4)\}$ and $\{(7,5),(3,3)\}$.

\begin{definition}
The number of alignements of a permutation $\sigma$ is equal to
$$
A_+(\sigma)+A_-(\sigma)+A_{+,-}(\sigma).
$$
\end{definition}
This definition looks a bit different but is equivalent to
the definition in \cite{wi}. The number of weak exceedances of a permutation
$\sigma$ is the cardinal of the set $\{j\ | \ \sigma(j)\ge j\}$. 
In \cite{wi} the following result was proved~:
\begin{proposition}
The number of permutations of $[n]$ with $k$ weak exceedances
and $\ell$ alignments is the coefficient of
$[q^{(k-1)(n-k)-\ell}]$ in 
$$
\hat{E}_{k,n}(q)=q^{k-k^2}\sum_{i=0}^{k-1}(-1)^i [k-i]_q^n q^{k(i-1)}\left( {n\choose i}q^{k-i}+{n\choose i-1}\right). 
$$
where $[k]_q=1+q+\ldots +q^{n-1}$.
\end{proposition}
These numbers have the property that if $q=1$ they are the eulerian numbers,
if $q=0$ they are the Narayana numbers and if $q=-1$ they are the binomial
coefficients. See \cite{wi} for details.

Lauren Williams proved in \cite{wi}
that ${E}(q,x,y)=\sum_{n,k}q^{n-k}\hat{E}_{k,n}(q)y^k x^n$ is equal to
$$
\sum_{i=0}^{\infty} \frac{y^i (q^{2i+1} - y)}
{q^{i^2+i+1} (q^i - q^{i+1}[i]x + [i]xy)}.
$$                                                                                
Here we exhibit the generating function  
$\hat{E}(q,x,y)=\sum_{n,k}\hat{E}_{k,n}(q)y^k x^n$
in a continued fraction form.
\begin{theorem}
$$
\hat{E}(q,x,y)=
\cfrac{1}{1-b_0x-\cfrac{\lambda_1x^2}{1-b_1x-\cfrac{\lambda_2x^2}
{1-b_2x-\cfrac{\lambda_3x^2}{\ldots}}}}
$$
with $b_n=y[n+1]_q+[n]_q$ and $\lambda_n=y[n]_q^2$
and $[n]_q=1+q+\ldots +q^{n-1}$.
\label{big}
\end{theorem}

To prove this theorem we introduce the notion of crossings. Let
\begin{itemize}
\item $C_+(i)=\{j\ | \ j<i\le\sigma(j)<\sigma(i)\}$; $C_+(\sigma)=\sum_{i=1}^n |C_+(i)|$.
\item $C_-(i)=\{j\ | \ j>i>\sigma(j)>\sigma(i)\}$, $C_-(\sigma)=\sum_{i=1}^n |C_-(i)|$.
\end{itemize}
We define the number of crossings of a permutation $\sigma$
to be $C_-(\sigma)+C_+(\sigma)$.

We can again define those parameters using the permutation diagram. $C_+(\sigma)$ is the number of pairs of
edges above the line that intersect and  $C_-(\sigma)$ is the number of pairs of
edges under the line that cross.

For $\sigma=(4,7,3,6,2,1,5)$ we have $C_{+}(\sigma)=2$ and $C_{-}(\sigma)=1.$
The pairs of edges
contributing to  $C_{+}(\sigma)$ are $\{(1,4),(4,6)\}$ and $\{(1,4),(2,7)\}$.
The pair of edges
contributing to  $C_{-}(\sigma)$ is $\{(7,5),(6,1)\}$.\\

The proof of the theorem is direct with the 
following two propositions.
\begin{proposition}
The coefficient $[x^ny^kq^\ell]$ in  $\hat{E}(q,x,y)$ is
the number of permutations of $[n]$ with $k$ weak exceedances 
with $\ell$ crossings.
\label{main1}
\end{proposition}

\begin{proposition}
For any permutation with $k$ weak exceedances
the number of crossings plus the number of alignments is
$(k-1)(n-k).$
\label{main}
\end{proposition}

We prove these two propositions in Section 2.
The proof of the first one uses a bijection of  Foata and Zeilberger 
between permutations and weighted bicolored Motzkin paths \cite{fz}.
The proof of the second one uses basic properties of permutations.
Then in Section 3 we use 
a bijection due to Fran\c{c}on and Viennot \cite{fv} 
and the bijection of Foata and Zeilberger and use 
results from \cite{cm,csz} to link with generalized patterns. 
One result is~:
\begin{proposition}
There is a one-to-one correspondance between permutations
with $k$ weak exceedances and $(n-k)(k-1)-\ell$ alignments 
and permutations with $k-1$ descents and $\ell$ occurrences
of the generalised pattern $13-2$.
\label{pattern}
\end{proposition}
The link between alignments and generalized
patterns was conjectured by Steingrimsson and Williams \cite{sw}.
In Section 4 we generalize to the case of decorated
permutations. In Section 5 we show how these numbers appear naturally
in the stationary distribution of the ASEP model.

\section{Proof of Proposition \ref{main1} and \ref{main}}

\subsection{Proof of Proposition \ref{main1}}

We use a bijection of Foata and Zeilberger \cite{fz} between permutations and weighted bicolored
Motzkin paths. We could also use
the bijection of Biane \cite{bi}. We refer to \cite{csz} for a compact definition
of these bijections.

A bicolored Motzkin path of length $n$ is a sequence
$c=(c_1,\ldots ,c_n)$ such that 
$c_i\in\{N,S,E,\bar{E}\}$ for $1\le i\le n$ such that
if 
$
h_i=\{ j<i\ |\ c_j=N\}-\{ j<i\ |\ c_j=S\}
$
then $h_i\ge 0$ for $1\le i\le n$ and $h_{n+1}=0$.

To any permutation $\sigma$ we associate a pair $(c,w)$ made of a bicolored Motzkin path
$c=(c_1,c_2,\ldots ,c_n)$ and a weight $w=(w_1,\ldots ,w_n)$.
The path is created using the following rules~:
\begin{itemize}
\item $c_i=N$ if $i<\sigma(i)$ and $i<\sigma^{-1}(i)$
\item $c_i=E$ if $i<\sigma(i)$ and $i>\sigma^{-1}(i)$
\item $c_i=\bar{E}$ if $i>\sigma(i)$ and $i<\sigma^{-1}(i)$
\item $c_i=S$ if $i>\sigma(i)$ and $i>\sigma^{-1}(i)$
\end{itemize}

The weight is created using the following rules~:
\begin{itemize}
\item $w_i=y p^{|A_+(i)|}q^{|C_+(i)|}$ if $c_i=N,E$.
\item $w_i=p^{|A_-(i)|}q^{|C_-(i)|}$ if $c_i=S,\bar{E}$.
\end{itemize}

For example $\sigma=(4,1,5,6,2,3)$ gives the path $(N,\bar{E},N,E,S,S)$
and the weight $(y,1,yp,yp,p,1)$.\\

We now need
\begin{lemma}
If $i\le \sigma(i)$ then
$$
|C_+(i)|=h_i-|A_+(i)|;
$$
and if $i>\sigma(i)$ then
$$
|C_-(i)|=h_i-1-|A_-(i)|.
$$
\end{lemma}
\noindent{\bf Proof.} It is easy to prove by induction that
$$
h_i=|\{j<i\ |\ \sigma(j)\ge i\}|=|\{j\ge i\ |\ \sigma(j)<i\}|.
$$
Using the definitions given in Section 1, 
if $i\le \sigma(i)$ then
$$
A_+(i)\cup C_+(i)=\{j<i\ |\ \sigma(j)\ge i\}
$$
and if $i>\sigma(i)$ then
$$
A_-(i)\cup C_-(i)=\{j>i\ |\ \sigma(j)<i\}=\{j\ge i\ |\ \sigma(j)<i\}\backslash\{i\}.
$$
This implies the result. \qed

Let ${\mathcal P}_n$ be the set of pairs $(c,w)$ obtained from permutations
of $[n]$.
Using the machinery developed in \cite{f,vi}, we get directly that if~:
$$
\hat{E}(q,p,x,y)=\sum_n x^n \sum_{(c,w)\in {\mathcal P}_n} \prod_{i=1}^n w_i.
$$
then
\begin{equation}
\hat{E}(q,p,x,y)=
\cfrac{1}{1-b_0x-\cfrac{\lambda_1x^2}{1-b_1x-\cfrac{\lambda_2x^2}
{1-b_2x-\cfrac{\lambda_3x^2}{\ldots}}}}
\label{best}
\end{equation}
with $b_n=y[n+1]_{p,q}+[n]_{p,q}$ and $\lambda_n=y[n]_{p,q}^2$.\\

This gives a generalization of Proposition \ref{main1}.
We define the number of nestings of a permutation $\sigma$ to be
$A_+(\sigma)+A_-(\sigma)$. 
\begin{proposition}
The coefficient of $[x^n y^k q^\ell p^m]$ in
$\hat{E}(q,p,x,y)$
is the number of permutations $\sigma$ of $[n]$ with $k$ weak 
exceedances, $\ell$ crossings
and $m$ nestings.
\end{proposition}

Note that this bijection implies that
\begin{proposition}
The number of permutations with $k$ weak exceedances and
$\ell$ crossings and $m$ nestings is equal to the number of permutations with $k$ weak exceedances and
$\ell$ nestings and $m$ crossings.
\end{proposition}

Similar results are known for set partitions and matchings \cite{kl,st}.

\subsection{Proof of Proposition \ref{main}}

We suppose that $\sigma$ is a permutation of $[n]$
with $k$ weak exceedances.
For any $i$ with $1\le i<n$, we first define~:
\begin{itemize}
\item $B_+(i)=\{j\ |\ j<i\le \sigma(j)\}$  
\item $B_-(i)=\{j\ |\ \sigma(j)< i\le j\}$  
\end{itemize}
Note that  $h_i=|B_+(i)|=|B_-(i)|$ and that for $i>\sigma(i)$, 
$
A_-(i)\cup C_-(i)\cup \{i\}=B_-(i).$ \\
Therefore
\begin{equation}
A_-(\sigma)+C_-(\sigma)=\sum_{i>\sigma(i)} (|B_-(i)|-1)=k-n+\sum_{i>\sigma(i)} |B_+(i)|
\label{eq1}
\end{equation}

For $i\le \sigma(i)$, let
$$
E_+(i)=\{j\ |\ i\in C_+(j)\}=\{j\ | \ i<j\le\sigma(i)<\sigma(j)\}.
$$ 
It is easy to see that~:
$$
E_+(i)\cup A_+(i)=B_+(\sigma(i)).
$$

Therefore
\begin{eqnarray*}
A_+(\sigma)+C_+(\sigma)&=& \sum_{i\le \sigma(i)}|{C}_+(i)|+ |A_+(i)| \\
&=& \sum_{i\le \sigma(i)}|E_+(i)|+ |A_+(i)| \\
&=&\sum_{i\le \sigma(i)}|B_+(\sigma(i))| \\
&=&\sum_{i>\sigma(i)} |D_+(i)| \label{eq2},
\end{eqnarray*}
where
$$
D_+(i)=\{j\ | \ j\le \sigma(j)\ {\rm and}\ i\in B_+(\sigma(j))\}.
$$
It is easy to see that for $i>\sigma(i)$, 
$D_+(i) =\{j\ |\ j\le \sigma(j)\ {\rm and}\ \sigma(i)<\sigma(j)<i\}$ 
and therefore that
$
B_+(i)\cup D_+(i)\cup A_{+,-}(i)=\{j\ |\ j\ge \sigma(j)\}.
$
As they are pairwise disjoint then
$$
|B_+(i)|+ |D_+(i)|+ |A_{+,-}(i)|=k. 
$$

Combining equations (\ref{eq1}) and (\ref{eq2})
we get that
$A_+(\sigma)+C_+(\sigma)+A_-(\sigma)+C_-(\sigma)+A_{+,-}(\sigma)=
k-n+ \sum_{i>\sigma(i)} |D_+(i)|+|B_+(i)|+|A_{+,-}(i)|.$
This concludes
the proof of Proposition \ref{main}. \qed

\section{Link with generalized patterns}

Continued fractions like the one presented in equation (\ref{best})
were studied combinatorially in \cite{cm,csz}.
A descent in a permutation is an index $i$ such that $\sigma(i)>\sigma(i+1)$.
An ascent in a permutation is an index $i$ such that $\sigma(i)<\sigma(i+1)$.
The pattern 31-2 (resp. 2-31, 13-2) occurs in $\sigma$ if there exist $i<j$ such that $\sigma(i)>\sigma(j)>\sigma(i+1)$ (resp.  $\sigma(j+1)<\sigma(i)<\sigma(j)$, $\sigma(i+1)>\sigma(j)>\sigma(i)$).\\

Theorem 10 in \cite{csz} associated with Theorem 22 in \cite{cm}
tells us that~:
\begin{proposition}
The coefficient of $[t^n x^k q^\ell p^m]$ in
$\hat{E}(t,x,q,p)$
is the number of permutations $\sigma$ of $[n]$ with $n-k$ descents, 
$\ell$ occurences of the patterns $31-2$ 
and $m$ occurences of the pattern $2-31$.
\end{proposition}
This can also be proved bijectively thanks to a bijection of Fran\c{c}on
and Viennot \cite{fv}. We present now that bijection. See also \cite{csz}.

Given a permutation $\sigma=(\sigma(1),\ldots , \sigma(n))$, we set $\sigma(0)=0$ and 
$\sigma({n+1})=n+1$. Let $\sigma(j)=i$.
Then $i$ is
\begin{itemize}
\item a valley if $\sigma({j-1})>\sigma({j})<\sigma({j+1})$
\item a double ascent  if $\sigma({j-1})<\sigma({j})<\sigma({j+1})$
\item a double descent  if $\sigma({j-1})>\sigma({j})>\sigma({j+1})$
\item a peak  if $\sigma({j-1})<\sigma({j})>\sigma({j+1})$
\item the beginning (resp. end) of a descent if $\sigma({j})>\sigma({j+1})$ 
(resp. $\sigma({j-1})>\sigma({j})$)
\item the beginning (resp. end) of a ascent if $\sigma({j})<\sigma({j+1})$ 
(resp. $\sigma({j-1})<\sigma({j})$).
\end{itemize}
If $\sigma(j)=i$ then we define $31-2(i)$ (resp. $2-31(i)$) to be the number of indices $k<j$ (resp. $k>j$)
such that $\sigma(k-1)>\sigma(j)>\sigma(k)$.\\

To any permutation, we associate a pair $(c,w)$ made of a bicolored Motzkin path
$c=(c_1,c_2,\ldots ,c_n)$ and a weight $w=(w_1,\ldots ,w_n)$.
The path is created using the following rules~:
\begin{itemize}
\item $c_i=N$ if $i$ is a valley
\item $c_i=E$ if $i$ is a double ascent
\item $c_i=\bar{E}$ if $i$ is a double descent
\item $c_i=S$ if $i$ is a peak
\end{itemize}

The weight is created using the following rules~:
\begin{itemize}
\item $w_i=y p^{31-2(i)}q^{2-31(i)}$ if $c_i=N,E$.
\item $w_i= p^{31-2(i)}q^{2-31(i)}$ if $c_i=S,\bar{E}$.
\end{itemize}

For example $\sigma=(6,2,1,5,3,4)$ gives the path $(N,\bar{E},N,E,S,S)$
and the weight $(y,1,yp,yp,y,1)$.\\

Now we prove the following Lemma
\begin{lemma}
For any $i$
$$
31-2(i)+2-31(i)=\left\{ \begin{array}{ll}
h_i & {\rm if}\ i \ {\rm is}\ {\rm the}\ {\rm beginning}\ {\rm of }\ {\rm an}\ {\rm ascent}\\ 
h_i-1 & {\rm if}\ i \ {\rm is}\ {\rm the}\ {\rm beginning}\ {\rm of }\ {\rm a}\  {\rm descent}\\ 
\end{array}\right.
$$
\end{lemma}
\noindent{\bf Proof.} 
We prove this lemma by induction. If $i$ is equal to 1 then 
$i$ is the beginning of an ascent and $31-2(1)+2-31(1)=0=h_1$.
 If $i>1$ then  $31-2(i)+2-31(i)=31-2(1)+2-31(i-1)+v$ where $v$ is zero,
 one or minus one.
It is easy to see that $v$ is one if $i-1$ is the end of a descent
and $i$ is the beginning of an ascent, $v$ is minus one if $i-1$ is the end of an ascent
and $i$ is the beginning of a descent and 0 otherwise. That gives exactly the lemma.  \qed\\

Let ${\mathcal P}_n$ be the set of pairs $(c,w)$ obtained from permutations
of $[n]$.
Using the machinery developed in \cite{f,vi}, we get directly that 
$\sum_n x^n \sum_{(c,w)\in {\mathcal P}_n} \prod_{i=1}^n w_i.$
is $\hat{E}(q,p,x,y)$ defined in equation (\ref{best}).\\

Combining the bijection of Fran\c{c}on and Viennot and the inverse of bijection
of Foata and Zeilberger \cite{fz} we get that 
\begin{proposition}
The number of permutations of $[n]$ with $k$ descents, 
$\ell$ occurences of the patterns $31-2$ 
and $m$ occurences of the pattern $2-31$
is equal to number of permutations $\sigma$ of 
$[n]$ with $n-k$ weak exceedances and $\ell$ crossings
and $m$ nestings.
\end{proposition}
 
We also propose the direct mapping. Starting from 
a permutation $\sigma$ with $k$ descents, 
$\ell$ occurences of the patterns $31-2$ 
and $m$ occurences of the pattern $2-31$, we form a permutation
$\tau$ with $k$ descents, 
$\ell$ occurences of the patterns $31-2$ 
and $m$ occurences of the pattern $2-31$. 

We first form two two-rowed arrays $f$ and $g$.
The first line of $f$ contains all the entries of $\sigma$
that are the beginning of a descent. They are sorted in increasing order. 
The second line of $f$ contains all the entries of $\sigma$
that are the end of a descent. They are sorted such that
if $i$ is in that row then $2-31(i)$ in $\sigma$ is the number
of entries of the row to the right of $i$ smaller than $i$.
The second line of $f$ contains all the entries of $\sigma$
that are the end of a descent sorted in increasing order.\\
The first line of $g$ contains all the entries of $\sigma$
that are the beginning of an ascent and that are sorted in increasing order. 
The second line of $g$ contains all the entries of $\sigma$
that are not the end of a descent. They are sorted such that
if $i$ is in that row then $2-31(i)$ in $\sigma$ is the number
of entries of the row to the left of $i$ greater than $i$.
We create the permutation $\tau$ which is the union of $f$ and $g$.\\

For example if $\sigma=(5,1,7,4,3,6,8,2)$ then the 
$2-31$ sequence is $$(2-31(1),\ldots ,2-31(8))=(0,0,1,1,2,1,1,0).$$
Then
$$
f=\left(\begin{array}{l}
4,5,7,8\\
1,3,4,2 \end{array}\right)
$$
and 
$$
g=\left(\begin{array}{l}
1,2,3,6\\
8,6,5,7
 \end{array}\right)
$$
Then $\tau=\left(\begin{array}{l}
1,2,3,4,5,6,7,8\\
8,6,5,1,3,7,4,2 \end{array}\right)
=(8,6,5,1,3,7,4,2)$.\\
 
We conclude this Section  by proving  Proposition \ref{pattern}.
Given a permutation $\sigma=(\sigma(1),\ldots ,\sigma(n))$.
Let $\pi=(\sigma(n),\ldots ,\sigma(1))$. If $\sigma$ has $k-1$ ascents
(or $n-k$ descents) and $\ell$ occurences of the pattern $31-2$, then $\pi$
has $k-1$ descents and $\ell$ occurences of the pattern $2-13$.
With the previous proposition this gives Proposition \ref{pattern}.

\section{Generalization for decorated permutations}

We can also derive the generating function of
$$
A_{k,n}(q)=\sum_{i=0}^{k-1}{n\choose i}E_{k,n-i}.
$$
These were introduced in \cite{wi}. 

The following corollary is an easy consequences of Theorem
\ref{big}.
Let $A(q,x,y)=\sum_{n,k}A_{k,n}(q)x^n y^k$.
\begin{corollary}
$$
A(q,x,y)=
\cfrac{1}{1-b_0x-\cfrac{\lambda_1x^2}{1-b_1x-\cfrac{\lambda_2x^2}
{1-b_2x-\cfrac{\lambda_3x^2}{\ldots}}}}
$$
with $b_n=(1+y)[n+1]_q$ and $\lambda_n=yq[n]_q^2$.
\end{corollary}

Lauren Williams  \cite{sw} proved that 
$$
A(q,x,y)=\frac{-y}{1-q}+\sum_{i\ge 1}
\frac{y^i(q^{2i+1}-y)}{q^{i^2+i+1}(q^i-q^i[i+1]x+[i]xy)}
$$

It would be interesting to have a direct proof of the identity
of the formal series form and the continued fraction form
of these generating functions \cite{sw}. \\

We can also interpret these results combinatorially.
In \cite{wi} the coefficient of $[q^{(n-k)k-\ell}]$ in $A_{k,n}(q)$
is interpreted in terms of decorated permutations with $k$
weak exceedances and $\ell$ alignments . Let us define these notions.
Decorated permutations
are permutations where the fixed points are bicolored \cite{po}.
We color these fixed points by colors $\{+,-\}$.
We say that $i\le_+ \sigma(i)$ if $i<\sigma(i)$ or $i=\sigma(i)$
and colored with color $+$. We say that 
$i\ge_- \sigma(i)$ if $i>\sigma(i)$ or $i=\sigma(i)$
and  colored with color $-$. 

For a decorated permutation $\sigma$ and $i$, let
\begin{itemize}
\item $A_+(i)=\{j\ |\ j<i\le_+ \sigma(i)<\sigma(j)\}$
\item $A_-(i)=\{j\ | \ j>i\ge_- \sigma(i)>\sigma(j)\}$
\item 
$A_{+,-}(i)=\{j\ |\ i\ge_- \sigma(i)>\sigma(j)\ge j\}\cup \{j\ |\ 
\sigma(j)\ge j>i\ge_- \sigma(i)\}$
\item $C_+(i)=\{j\ | \ i<j\le\sigma(i)<\sigma(j)\}$
\item $C_-(i)=\{j\ | \ j>i>\sigma(j)>\sigma(i)\}$
\end{itemize}
As for permutations we define ${\rm A}_{\pm}(\sigma)=\sum_i {\rm A}_{\pm}(i)$.

With these notions, we can again define the number of alignments (resp. nestings, crossings)
of a decorated permutation $\sigma$ as $A_+(\sigma)+A_-(\sigma)+
A_{+,-}(\sigma)$ (resp. $A_+(\sigma)+A_-(\sigma)$, $C_+(\sigma)+C_-(\sigma)$). 
The number of weak exceedances (resp. descedances) of a decorated permutation
is the cardinality of the set $\{i\ |\ i\ge_+ \sigma(i)\}$ (resp. $\{i\ |\ i< \sigma(i)\}$ ).

Let
$$
A(q,p,x,y)=\cfrac{1}{1-b_0x-\cfrac{\lambda_1x^2}{1-b_1x-\cfrac{\lambda_2x^2}
{1-b_2x-\cfrac{\lambda_3x^2}{\ldots}}}}
$$
with $b_n=(1+y)[n+1]_{p,q}$ and $\lambda_n=yq[n]_{p,q}^2$.

A direct generalization of the bijection \`a la 
Foata-Zeilberger on decorated permutations gives~:
\begin{proposition}
The coefficient of $[x^n y^k q^\ell p^m]$ in
$A(q,p,x,y)$
is the number of decorated permutations $\sigma$ of $[n]$ with $k$ weak 
exceedances, $\ell$ is the sum of the crossings and the descedances
and $m$ nestings.
\end{proposition}

We can also make a direct link with the alignments~:
\begin{proposition}
For any decorated permutation $\sigma$ with $k$ weak exceedances
$$
A_+(\sigma)+A_-(\sigma)+C_+(\sigma)+C_-(\sigma)+A_{+,-}(\sigma)
+|\{j\ |\ j>\sigma(j)\}|=(n-k)k.
$$
\end{proposition}
\noindent{\bf Proof.} The proof is omitted as it follows exactly
the same steps as the proof of Proposition \ref{main}. \qed \\

We leave as an open question to derive the bijection \`a la 
Fran\c{c}on-Viennot on decorated permutations to interpret the
$A_{k,n}(q)$ in terms of descents and generalized patterns.

\section{Link with the ASEP model}

The ASEP model \cite{mg,der} 
consists of black particles entering a row of $n$
cells, each of which is occupied by a black particle or vacant. A
particle may enter the system from the left hand side, hop to the
right or to the left and leave the system from the right hand side, with the
constraint that a cell contains at most one particle. 
We will say that the empty cells are filled with white
particles $\whi$. A basic configuration is a row of $n$ cells, each
containing either a black $\bla$ or a white $\whi$ particle. 
Let $\basic_n$ be the
set of basic configurations of $n$ particles. We write these
configurations as though they are words of length $n$ in the language
$\{\whi,\bla\}^*$.

The ASEP defines a Markov chain $P$ defined on $\basic_n$ with the
transition probabilities $\alpha$, $\beta$, 
and $q$.  The probability $P_{X,Y}$, of finding the system in state
$Y$ at time $t+1$ given that the system is in state $X$ at time $t$ is
defined by:
\begin{subequations}
\begin{itemize}
\item If $X=A\bla\whi B$ and $Y=A\whi\bla B$ then
\begin{equation}
P_{X,Y}=1/(n+1);\ \ \ \
P_{Y,X}=q/(n+1)
\end{equation}
\item If $X=\whi B$ and $Y=\bla B$ then
\begin{equation}
P_{X,Y}=\alpha/(n+1).
\end{equation}
\item If $X=B\bla$ and $Y= B\whi $ then
\begin{equation}
P_{X,Y}=\beta/(n+1).
\end{equation}
\item Otherwise $P_{X,Y}=0$ for $Y\neq X$ and $P_{X,X}=1-\sum_{X\neq
    Y}P_{X,Y}$.
\end{itemize}
\end{subequations}

\begin{figure}[ht!]
  \begin{center}
    \includegraphics[scale=0.5]{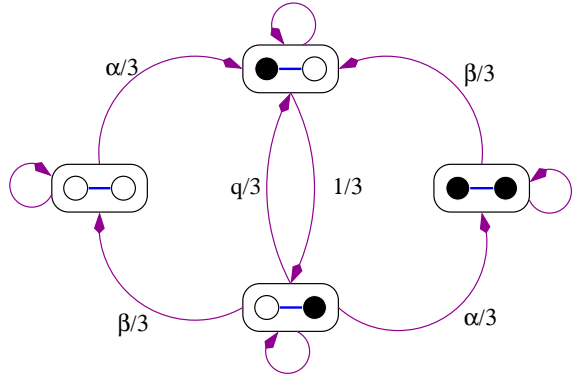}
  \end{center}\centering
  \caption{The chain $P$ for $n=2$.}
  \label{chainp}
\end{figure}

See an example for $n=2$ in Figure~\ref{chainp}.

This Markov chain has a unique stationary distribution \cite{der}.
The case $q=0$ was studied combinatorially by Enrica Duchi and Gilles Schaeffer
\cite{ds}. In \cite{bcer} another combinatorial approach was taken
to treat the general case.

\begin{definition}\cite{bcer}
  Let $\mathcal{P}(n)$ be the set of bicolored Motzkin paths of
  length $n$. The weight of the path in $\mathcal{P}(0)$ is 1.  The
  weight of any path $p$ denoted by $w(p)$ is the product of the weights of
  its steps. The weight of a step, $p_i$, starting at height $h$ is
  given by:
  \begin{subequations}
  \begin{eqnarray*}
    \mbox{if } p_i = N & \mbox{ then } & w(p_i) = [h+1]_q\\
    \mbox{if } p_i = \bar{E} & \mbox{ then } & w(p_i) =[h]_q+q^h/\alpha \\
    \mbox{if } p_i = E & \mbox{ then } & w(p_i) =[h]_q+q^h/\beta \\
    \mbox{if } p_i = S & \mbox{ then } & w(p_i) =[h]_q
+q^{h}/(\alpha\beta)-q^{h-1}(1/\alpha-1)(1/\beta-1).
  \end{eqnarray*}
\end{subequations}
\end{definition}

Given a path $p$, $\theta(p)$ is the basic configuration
such that each $\bar{E}$ and $S$ step is changed to $\whi$ and each
$E$ and $N$ step is changed to $\bla$. Let
\begin{equation}
  W(X) = \sum_{p \in \theta^{-1}(X)} w(p)
\end{equation}
and 
\begin{equation}
  Z_n = \sum_{X\in B_n} W(X)
\end{equation}

\begin{theorem}\cite{bcer}
At the steady state, the probability that the chain is in the basic 
configuration $X$ is
$$
\frac{W(X)}{Z_n}
$$
\end{theorem}

We can use these results and observations from the previous
sections to get~:
\begin{theorem}
If $\alpha=\beta=1$, at the steady state, the probability that 
the chain is in a basic configuration with $k$ particles is
$$
\frac{\hat{E}_{k+1,n+1}(q)}{Z_n}.
$$
\end{theorem}

Before proving that theorem, we need a Lemma
\begin{lemma}
There is a weight preserving bijection between $\mathcal{P}({n,k})$ the set of bicolored 
Motzkin paths of length $n$ where the  weight  of any step starting at height $h$ is $[h+1]_q$
and where $k$ is the number of steps $N$ plus the number of steps
$E$ 
and  $\mathcal{P}'({n+1,k+1})$ the set of Motzkin paths of length $n+1$
where the weight of any step starting at height $h$ is
 $[h+1]_q$ if the step is $N$ or $E$ and  $[h]_q$ otherwise
 and where $k+1$ is the number of steps $N$ plus the number of steps
$E$.  
\end{lemma}
\noindent{\bf Proof.} This is  a general version of the classical bijection
between bicolored Motzkin path of length $n$ and bicolored Motzkin path
of length $n+1$ where one of the horizontal steps can not appear at height 0.
The proof consists of transforming each step 
$N$ (resp. $E$, $\bar{E}$, $S$) weighted by $w$ 
by two steps $(N,N)$ (resp. $(N,S)$, $(S,N)$, $(S,S)$)
each of them weighted by $w/2$. 
Then we add a $N$ step at the beginning and add a step $S$ at the end 
both with weight $1/2$ and apply the reverse map. \qed\\

\noindent{\bf Proof of the Theorem.} 
We want to compute for $\alpha=\beta=1$
$$
W({k,n})=\sum_{\begin{subarray}{c}X\in B_{n}\\
X \ {\rm has} \ k\ {\rm particles}\end{subarray}}W(X)=\sum_{p\in \mathcal{P}(n,k)}w(p).
$$
Now we use the previous lemma and get
$$
W({k,n})=\sum_{p\in \mathcal{P}'(n+1,k+1)}w(p).
$$
Using \cite{f,vi} and the definition of the weight of the steps of $\mathcal{P}'(n+1,k+1)$,
we  conclude that~:
$$
1+\sum_{n\ge 0} x^{n+1} \sum_{k=0}^{n} y^{k+1} W({k,n})=\hat{E}(q,x,y)
$$
and therefore that  $W({k,n})=\hat{E}_{k+1,n+1}(q)$.
\qed\\

\section{Conclusion}

Several open problems naturally arise~:
\begin{itemize}
\item Can we define generalized patterns for decorated permutations?
\item Can we generalize these $q$-Eulerian numbers to understand the ASEP
when $\alpha\neq 1$ or $\beta\neq 1$? 
\item Can we extend the definition of $k$-crossing and $k$-nesting that were
defined for matchings and set partitions \cite{st} ?
\end{itemize}

\noindent{\bf Acknowledgments.} The author wants to thank Lauren Willians and
Petter Br\"and\'en for their encouragements and constructive comments and
the Institute Mittag Leffler where this work was done.

\end{document}